# A Method for Generating a Well-Distributed Pareto Set in Nonlinear Multiobjective Optimization

S.V. Utyuzhnikov[1,2], P. Fantini[1], and M.D. Guenov[1]

[1]*Department PPAE, School of Engineering Cranfield University,*
*Cranfield, Bedfordshire, MK43 0AL, UK*
*e-mail: m.d.guenov@cranfield.ac.uk*

[2]*School of Mechanical, Aerospace & Civil Engineering,*
*University of Manchester, Sackville Street, M60 1QD*
*e-mail: s.utyuzhnikov@manchester.ac.uk*

In multidisciplinary optimization the designer needs to find solution to optimization problems which include a number of usually contradicting criteria. Such a problem is mathematically related to the field of nonlinear vector optimization where there are many numerical methods capable of providing a solution. However, only a few of those are suitable for real multidisciplinary design in industry because an iteration design circle usually is very time-consuming. This is due to the time scales and computational resources associated with each iterative design cycle. The recently suggested Physical Programming Method appears to match many requirements raised in industrial applications. The method is modified to make its realization easier and more efficient, the main focus being the even generation of the complete Pareto set. The method is used to find the Pareto surface for different test cases.

## Introduction

In real industrial design, the decision-maker (DM) has to take into account many different criteria such as low initial cost, manufacturability, long life and good performance which cannot be satisfied simultaneously. In fact, it is possible only to consider a trade-off among all (or almost all) criteria. The task becomes even more complicated because of additional constraints which always exist in practice.

Mathematically, the trade-off analysis can be formulated as a vector nonlinear optimization problem under constraints. Generally speaking, the solution of such a problem is not unique. It is naturally to exclude from the consideration any design solution which can be improved without deterioration of any discipline and violation of the constraints; in other words, a solution which can be improved without any trade-off. This leads to the Pareto optimal solutions[1].

Mathematically, each Pareto point is a solution of the multidisciplinary optimization problem. A designer selects the ultimate solution among the Pareto set on the basis of additional requirements (which may be subjective). Generally speaking, it is desirable to have a sufficient number of Pareto points to represent the whole Pareto set. Yet, it is important that the Pareto set is to be evenly distributed, otherwise the representation of the Pareto surface becomes inefficient.

In spite of the existence of many numerical methods for vector optimization, there are few methods suitable for real-design industrial applications, especially for preliminary design. This is due to high time-limit requirements. In many practical MDO applications the design cycle includes time-consuming and expensive computations at each discipline. In the aerospace industry this is most evident in the solution of the aerodynamics and stress analysis tasks. The solutions corresponding to these subtasks influence each other and usually demand iterations to reach the general solution. Because of the different mathematical properties of the governing equations describing these disciplines and their relatively weak coupling, it is not efficient to solve these equations simultaneously. Under such conditions, it is necessary to minimize the number of iterations required to find a Pareto optimal solution by the choice of an appropriate numerical optimization method. Also, it is very desirable to utilise the experience of a DM when considering issues which are beyond formalization, but can be very important from practical point of view.

Usually, vector optimization methods lead to a scalar optimization of an aggregated objective function (AOF) which includes a combination of objective (cost) functions. The most distributed case of the AOF is represented by a linear (weight) combination of the objective functions[1]. Although this method is relatively simple, it has many drawbacks related mainly with the uncertainty about weights. It may require many iterations to find the combination of the weights leading to a solution which corresponds to the DM's expectations[2]. Furthermore, it is well known that this method can generate only the convex part of a Pareto surface[3] while structural problems often result in non-convex Pareto surfaces. One may note that this drawback can be avoided by using either a more complex consideration of the AOF[4] or the weighted Tchebycheff method[1]. Nevertheless, the weights remain to be unknown functions of the objectives[2].

In real industrial design, especially in preliminary design, the DM is able to consider only a few possible solutions (Pareto points). In such a context, it is important to have an even distribution of Pareto points to obtain maximum information on the Pareto surface at minimum (computational) cost. Das and Dennis[5] showed that an even spread of weights in the AOF does not necessary result in an even distribution of points in the Pareto set. Also, the spread of the points strongly depends on the relative scaling of the objectives. In[6,7,8], necessary conditions for an AOF were obtained for capturing any Pareto point. In the authors' knowledge, there are only three methods which are able to guarantee an even distribution of a whole Pareto surface[9].

The Normally Boundary Intersection (NBI) Method[10,11] was developed for generating an even distribution of Pareto points by Das and Dennis. The method appeared to generate non-Pareto and locally Pareto points that requires a filtering procedure[12,13]. The new Normal Constraint (NC) Method[9,12] developed recently looks very promising. Both methods have clear geometrical interpretation. They are based on the well-known fact that a Pareto surface belongs to the boundary of the feasible space towards minimization of the objective functions[1]. So-called *anchor points* are obtained in the feasible objective space, first. An anchor point corresponds to the optimal value of one and only one objective function in the feasible space. Thus, $n$ objective functions give $n$ anchor points. Second, the *utopia plane* passing through the anchor points is obtained. In both the NBI and NC methods, the Pareto surface is then obtained by the intersection of lines normal to the utopia plane and the boundary of the feasible space. The single optimization problem, used in the NC, is based only on inequality constraints. This modification makes the method more flexible and stable. An even distribution of Pareto points is provided by even distribution of the lines orthogonal to the utopia plane. Although, one can note that this is valid until cosine of the angle between such a line and the normal to the boundary of the feasible space is not locally close to zero. Both methods may fail to generate Pareto solutions over the entire Pareto frontier[9] in multidimensional case. The recent modification of the NC[9] eliminates this drawback and guarantees the complete representation of a Pareto frontier. Meanwhile, both methods may generate non-Pareto and locally Pareto solutions[12] though the NC does it less likely[9].

The Physical Programming (PP) Method was suggested by Messac in[14]. This method also generates Pareto points on both convex and non-convex Pareto frontiers as it was shown in[15,16]. The method does not use any weight coefficients and allows one to take into account the DM experience immediately. In this sense, it looks to be the most interesting method for practical applications under the above stated conditions. In the PP, the designer assigns each objective to one of the four categories (class-functions). The optimization is based on minimization of an aggregate preference function determined by the preference functions (class-functions) with preferences set *a priori*. The notion of the generalized Pareto optimal solution was introduced in[17] on the basis of the class-functions.

The PP is modified below to make it simpler and more efficient for practical applications. A simpler structure of the class-functions is derived. The class-functions are generalized to shrink the search domain and make its location in space more optimal. This is critical when generating an even set of the Pareto frontier. The proposed modification allows us to combine the advantages of the PP, NBI and NC methods. The algorithm to obtain an even distribution of the Pareto set is described in this paper. One of its main advantages is that it does not provide non-Pareto solutions while local Pareto solutions may be easily recognized and removed. A comparison against the standard approach shows the modified approach is able to generate a much more even Pareto set. The suggested method is then generalized for the search of the generalized Pareto solutions.

**I. Multiobjective optimization problem. Pareto optimization**

It is assumed that an optimization problem is described in terms of a design variable vector $\mathbf{x} = (x_1, x_2, ..., x_N)^T$ in the design space $\mathbf{X} \in R^N$. A function $\mathbf{f} \in R^M$ evaluates the quality of a solution by assigning it to an objective vector $\mathbf{y} = (y_1, y_2, ..., y_M)^T$ ($y_i = f_i(\mathbf{x})$, $f_i: R^M \to R^1$, $i = 1, 2, ..., M$) in the objective space $\mathbf{Y} \in R^M$. Thus, $X$ is mapped by $f$ onto $Y$: $X \mapsto Y$. A multiobjective optimization problem may be formulated in the following form:

$$\text{Minimize } [y(x)] \quad (1)$$

subject to $K$ inequality constraints

$$g_i(x) \leq 0, \quad i = 1, 2, ..., K \quad (2)$$

and $P$ equality constraints

$$h_j(x) = 0, \quad j = 1, 2, ..., P \quad (3)$$

The feasible design space $X^*$ is defined as the set $\{x / g_j(x) \leq 0, j = 1, 2, ..., K;$ and $h_i(x) = 0, i = 1, 2, ..., P\}$. The feasible criterion (objective) space $Y^*$ is defined as the set $\{Y(x)/ \mathbf{X} \in \mathbf{X}^*\}$. The feasibility means no constraint is violated.

A design vector $\mathbf{a}$ ($\mathbf{a} \in \mathbf{X}^*$) is called a Pareto optimum iff it does not exist any $\mathbf{b} \in \mathbf{X}^*$ such that

$$\mathbf{y}(\mathbf{b}) \leq \mathbf{y}(\mathbf{a}) \text{ and exist } l \leq M: y_l(\mathbf{b}) < y_l(\mathbf{a}).$$

A design vector is called a local Pareto optimum if it is a Pareto optimum within its some neighbourhood.

**II. Physical programming method**

In the PP, both objective and constraint functions are treated as design metrics. Each design metric $F_i(x)$ is refereed to one of possible (determined in advance) utility dimensionless functions $\bar{F}_i(F_i(\mathbf{x}))$ called *class functions*. The determination of the class-functions reflects the qualitative classifications of possible preferences. The preferences are split into two major groups: hard class-functions *H* and soft ones *S*. The former reflect the presence of constraints (2) and (3), while the latter ones reflect the preferences settled for objective functions which may include subjective preferences of the DM based, say, on his/her experience and intuition.

There are considered four soft class-functions: the *1S* (smaller is better) where an objective is to be minimized; the mirror function *2S* (larger is better) where an objective is to be maximized; the class *3S* (value is better) where an objective is to be close to a particular preferred value; and the class *4S* (range is better) which is similar to the *3S* function but the preferable value lies in some range rather than corresponds to some value selected *a priori*. All soft class-functions are positive and dimensionless. The argument of a soft class-function is subdivided by the DM into different preference ranges: highly desirable, desirable, tolerate, undesirable, highly undesirable and unacceptable. Such a subdivision is a part of the approximation of the class-functions and allows the DM to exploit own experience. The values of the class-functions at the boundaries of the ranges are fixed. Therefore, scaling between different objective functions is automatically provided.

The qualitative behaviour of the *S* class-functions is given in Figure 1. The class-functions *1S* and *2S* tend to zero if the argument tends to minus or plus infinity, accordingly. Each of the class-functions *3S* and *4S* has one and only one minimum point. All class-functions have the same values at the boundary points of the preference ranges regardless of the type of a design metric. Only the locations of the ranges depend on a metric. The presence of the class-functions *3S* and *4S* does not violate the generality of formulation (2)-(4). The optimization problem can be reformulated as minimization problem (2)-(4). Yet, the PP lexicon makes the formulation closer to the formulation of real design optimisation problem. In fact, the whole approach is more reflective of real life[14]. In particular, the classes *3S* and *4S* represent very common in practice cases when it is impossible to determine if a smaller or larger design metric is better or worse.

As it is easy to see from Figure 1, the more desirable range the less value of the class-function. The preference between the ranges is enforced, that means: the highly undesirable region significantly worse that the undesirable, the undesirable region significantly worse than the tolerable region, and so on. This demand is reflected by the following conditions[14]:

$$\Delta_k \overline{F}_i > n_{sc} \Delta_{k-1} \overline{F}_i, \qquad 1 \leq k \leq 4, \qquad (4)$$

where $\Delta_k \overline{F}_i = \overline{F}_i(F_{i,k+1}) - \overline{F}_i(F_{ik})$, $\Delta_0 \overline{F}_i = \overline{F}_i(F_{i1})$, $F_{ik}$ $(k=1, 2, ..., 4)$ are the boundary points of the preference ranges for a metric $F_i$; $i$ is a number of a soft-class function, $n_{sc}$ is the number of soft design metrics.

Ultimately, the problem is reduced to the following optimization problem under constraints:

$$\min G(\mathbf{x}) = \frac{1}{n_{sc}} \log_{10}[\sum_{i=1}^{n_{sc}} \overline{F}(F_i)], \qquad (5)$$

subject to:

$$\begin{aligned}
F_i(\mathbf{x}) \leq F_{i5} & \qquad \text{(for } 1S\text{)} & (6) \\
F_i(\mathbf{x}) \geq F_{i5} & \qquad \text{(for } 2S\text{)} \\
F_{i5L} \leq F_i(\mathbf{x}) \leq F_{i5R} & \qquad \text{(for } 3S\text{)} \\
F_{i5L} \leq F_i(\mathbf{x}) \leq F_{i5R} & \qquad \text{(for } 4S\text{)}
\end{aligned}$$

and constraints related with the hard classes

$$\begin{aligned}
F_i(\mathbf{x}) \leq F_{iM} & \qquad \text{(for } 1H\text{)} & (7) \\
F_i(\mathbf{x}) \geq F_{im} & \qquad \text{(for } 2H\text{)} \\
F_{im} \leq F_i(\mathbf{x}) \leq F_{iM} & \qquad \text{(for } 3H\text{)}
\end{aligned}$$

The equality constraints (3) can be changed by double inequality constraints of *1H* and *2H* types:

$$h_j(\mathbf{x}) \leq 0,$$
$$h_j(\mathbf{x}) \geq 0$$

Any constraints on the design variables *x* are considered as ordinary constraints. The logarithm in (5) is only used to diminish the difference between the maximal and minimal values that may effect on the convergence of an iteration algorithm. It is easy to see that the more "undesirable" range is, the more it effects on the value of the AOF. Thus, preferences play somewhat the same role as weights in the weighted-sum method[14]. In (5) and henceforth, we assume that in the formulation $\overline{F}(F_i)$ the appropriate category of the class-function is used. This means that if an objective $F_i$ is to be minimized then *1S* class-function is utilised, if $F_i$ is to be maximized then *2S* class-function is considered and so on.

A robust implementation of the PP necessitates that problem (5) - (7) has a unique solution, which means absence of local minima. To guarantee this, each class-function must be a strictly convex function[14], that is:

$$\frac{d\overline{F}}{dF_i} \subset C^1, \quad \frac{d^2\overline{F}}{dF_i^2} > 0 \qquad (8)$$

It is easy to see that both high-order polynomial and cubic spline functions are not acceptable. Messac[14] derived the class-functions that satisfy all requirements mentioned above. They have been implemented into the PP and successfully used for solving many MDO problems[18,19]. We suggest another approximation of the class-functions which also satisfy all the requirements, but look simpler and more compact.

Let us assume that the derivative of the class-function $1S$ has the following form:

$$\frac{d\overline{F}}{dF} = Ae^{\alpha(F)}, \qquad (9)$$
$$A > 0,$$
$$\alpha'(F) > 0,$$

where index $i$ is omitted for the sake of simplicity. Then, function $\overline{F}$ is strictly convex because both the first and the second derivatives are positive.

Further, $\alpha$ must be a smooth function which satisfies condition (4). Since $\Delta_k \overline{F} = A \int_{F_k}^{F_{k+1}} e^{\alpha(F)} dF$, the latter condition is equivalent to

$$A\Delta_k F \int_0^1 e^{\alpha^{(k)}(\xi)} d\xi = \Delta_k \overline{F}, \qquad (10)$$

where $\alpha^{(k)} = \alpha^{(k)}(\xi^{(k)})$, $\xi^{(k)} = \dfrac{F - F_k}{F_{k+1} - F_k}$ $(F_k \le F \le F_{k+1})$.

Additional condition on $\alpha$ follows from (8), (9) and means:

$$\alpha^{(k)}(0) = \alpha^{(k-1)}(1) \qquad (11)$$

It follows that the linear function

$$\alpha^{(k)}(\xi) = a_k \xi + b_k \qquad (12)$$

is able to satisfy conditions (10) and (11) if

$$a_k = A\Delta_k F e^{b_k}(e^{a_k} - 1)/\Delta_k \overline{F}, \quad (k = 1, 2, ..., 4) \qquad (13)$$
$$b_k = a_{k-1} + b_{k-1}$$

Conditions (13) give a recurrent relation for the calculation of $a_k$ and $b_k$, where the former condition assumes solution of a nonlinear equation.

The values of $\overline{F}$ at the boundaries of the ranges do not depend on $F_{ik}$ and are fixed using relation

$$\Delta_k \overline{F} = \beta n_{sc} \Delta_{k-1} \overline{F}, \qquad (k = 1, 2, 3, 4) \qquad (14)$$
$$\beta > 1, \ \Delta_0 \overline{F} = 1$$

Having integrated (9), the ultimate solution is obtained as follows

$$\overline{F}^{(k)} = \overline{F}_k + \Delta_k \overline{F} \frac{e^{a_k \xi^{(k)}} - 1}{e^{a_k} - 1} \qquad (F_k \leq F \leq F_{k+1}), \qquad (15)$$

where $\xi^{(k)} = \dfrac{F - F_k}{F_{k+1} - F_k}$, $\overline{F}_k = \overline{F}(F_k)$, $k = 1, 2, 3, 4$.

In the rest region, $F < F_1$, it is possible to choose the following exponential function

$$\overline{F}^{(0)} = e^{A(F - F_1)} \qquad (16)$$

then the initial conditions in (13) may be chosen as uniform:

$$a_0 = 0, \ b_0 = 0. \qquad (17)$$

As we see, $\overline{F}^{(0)} = 1$ always.

The parameter $A$ must have the dimensionality of $F^{-1}$. It is naturally to set it as

$$A = \frac{1}{F_5 - F_1}. \qquad (18)$$

Thus, the class-function $\overline{F}$ is given by (15) and (16) where the coefficients are determined by (13), (17) and (18). Obviously, the function $\overline{F}$ is positive. From (16), it follows that $\overline{F}(F_1) = 1$.

Finally, we will give a brief comment on solving the nonlinear equation with respect to $a_k$ in (13). The equation can be rewritten as follows (index $i$ is omitted):

$$a_k = \varphi(a_k), \qquad (19)$$

where $\varphi(a_k) = \ln(1 + a_k / \omega_k)$, $\omega_k = A \Delta_k F e^{b_k} / \Delta_k \overline{F}$.

*Proposition 1*: A positive solution of (13) always exist and the method of simple iteration

$$a_k^{(n+1)} = \varphi(a_k^{(n)}) \qquad (n = 0, 1, 2, ....) \qquad (20)$$

monotonically converges to this solution at any initial $a_k^{(0)} > -\ln \omega_k$.

*Proof*:

From

$$\Delta_k \overline{F} = \int_{F_k}^{F_{k+1}} \frac{d\overline{F}}{dF} dF > \int_{F_k}^{F_{k+1}} \frac{d\overline{F}}{dF}\bigg|_k dF = Ae^{b_k}\Delta_k F,$$

a lower evaluation of $\omega_k$ follows:

$$\omega_k < 1.$$

Then $\dfrac{d\varphi}{da}(0) = \omega_k^{-1} > 1$ and, apart from 0, the second solution of (13) always exists, it being positive.

At the positive solution $a_k = a_k^*$ it is possible to obtain an upper evaluation of $\omega_k$:

$$\Delta_k \overline{F} = \int_{F_k}^{F_{k+1}} \frac{d\overline{F}}{dF} dF < \int_{F_k}^{F_{k+1}} \frac{d\overline{F}}{dF}\bigg|_{k+1} dF = Ae^{a_k+b_k}\Delta_k F,$$

so that:

$$\omega_k > e^{-a_k}.$$

Thus, the following evaluation is valid

$$0 < \frac{d\varphi}{da_k}(a_k^*) = \frac{1}{\omega_k e^{a_k^*}} < 1.$$

It is easy to see that at any point $a_k > a_k^*$: $\varphi'(a_k) < \varphi'(a_k^*)$. This means, mapping (19) appears to be contractive and iteration process (20) always monotonically converges if the initial approximation $a_k^{(0)} > -\ln \omega_k$ that immediately follows from the Banach fixed point theorem[20]. □

The class-function *2S* is obtained as the mirror-function of *1S*. The class-functions *3S* and *4S* can be created in regions $F < F_{1L}$ and $F > F_{1R}$ using the class-functions *2S* and *1S*, accordingly. In the rest region $F_{1L} < F < F_{1R}$, the class functions are defined as follows:

$$\overline{F}^{(0)}(F) = A_0 \xi^{(0)}(F)^{2m} + \varepsilon_0, \qquad (21)$$

where

$$\xi^{(0)}(F) = \frac{2F - F_{1L} - F_{1R}}{F_{1R} - F_{1L}},$$

$$A_0 = \frac{F_{1R} - F_{1L}}{4m} A, \quad \varepsilon_0 = 1 - A_0.$$

Integer parameter *m* defines the "flatness" of the function in the vicinity of the point of the minimum. It is naturally to set m >> 1 for 4S class-function. In the case of the 3S class-function, $F_{1L}$ should be close to $F_{1R}$, and $F_1 = (F_{1R} + F_{1R})/2$. Formula (21) then represents a local smoothness at point $F_1$. Then, $\overline{F}(F_1)$ approximately equals 1. Yet, it is not difficult to reach the exact equality for *3S* class-functions at $F_1$ using the following function in the interval $F_{1R} \leq F \leq F_{2R}$:

$$\overline{F}^{(1)} = A_1(e^{a_1[\xi^{(1)}]^2} - 1) + 1,$$

$$A_1 = \frac{\overline{F}_2 - 1}{e^{a_1} - 1}, \quad a_1 = \frac{(F_{2R} - F_{1R})A}{2A_1}.$$

In this case, in (13) $b_2 = a_1$. The similar solution can be obtained in the interval $F_{2L} \leq F \leq F_{1L}$.

### III. Generation of even distribution of a Pareto set

The PP allows one to generate an even distribution of the Pareto frontier. The appropriate algorithm is given in[16]. The original approach is briefly described below, followed by the description of our proposed modification which aims to make the algorithm more efficient, especially in the case of a concave Pareto frontier. Also, the algorithm is generalized on the *2S – 4S* class-functions.

Let us define the trade-off matrix *T* as follows:

$$T = \begin{bmatrix} F_{1,min} & F_{12} & \ldots & F_{1n_{sc}} \\ F_{21} & F_{2,min} & \ldots & F_{2n_{sc}} \\ \ldots & \ldots & \ldots & \ldots \\ F_{n_{sc}1} & F_{n_{sc}2} & \ldots & F_{n_{sc},min} \end{bmatrix}, \quad (22)$$

In the trade-off matrix *T*, an *i*-th row represents the coordinates of an anchor point $\boldsymbol{\mu}_i^*$ corresponding to the solution of single-optimization problem *min $F_i$* in the feasible criterion space $\mathbf{Y}^*$.

In the feasible space $\mathbf{Y}^*$ a hypercube *H* limiting the search domain is defined in the following manner. We set the pseudo nadir point[9] $F_{i,max} = \max_j F_{ij}$ that corresponds to the maximum *i*-th component among all anchor points. Then, the hypercube *H* is represented as follows: $H = [F_{1,min} F_{1,max}] \times [F_{2,min} F_{2,max}] \times \ldots \times [F_{n_{sc},min} F_{n_{sc},max}]$.

For the sake of simplicity, consider the 2D case where there are only two design metrics. We assume that each of the design metrics belongs to class *1S*. Following[16], for each of the design metrics let us introduce the vector of pseudo-preferences $P_i$:

$$P_i \equiv (F_{i1}, F_{i2}, F_{i3}, F_{i4}, F_{i5})^T = F_i^{(0)}(1,1,1,1,1)^T + a_i(0, \frac{1}{4}, \frac{1}{2}, \frac{3}{4}, 1)^T, \quad (23)$$

where $a_i = (F_{i,\max} - F_{i,\min})/n_d$, $F_i^{(0)}$ is a free parameter. The parameter $n_d$ defines the box size.

In such a case, either region $F_1 > F_{15}$ or region $F_2 > F_{25}$ becomes unacceptable. Thus, the box $D = [F_{11} F_{15}] \times [F_{21} F_{25}]$ defined by the pseudo-preferences limits the search domain from the right and upper sides leaving it in the other directions unlimited for a formal search. Changing the free vector $\mathbf{F}^{(0)} = (F_1^{(0)}, F_2^{(0)}, ..., F_{n_{sc}}^{(0)})^T$, it is possible to shift the box $D$ in the hypercube $H$ to seek Pareto solutions. The current location of the box $D$ determines a possible location of a Pareto point since the Pareto points outside $D$ (more precisely – higher or on the right of $D$) are excluded from the current consideration. In some sense, it operates similarly to the $\varepsilon$-constraint method[1], but in contrast, space reduction is simultaneously performed for all objectives. For example, moving the box to the lower-right angle of the hypercube $H$ (in 2D case), we give preferences to low values of the second objective at the expense of high values of the first objective. In[16], the algorithm is given for shifting the box $D$ over the space $\mathbf{Y}^*$ to seek the Pareto frontier. To provide this, the free vector $\mathbf{F}_i^{(0)}$ is specially chosen to move the box along lines parallel to a diagonal of the hypercube $H$ passing trough the lower-right angle and the upper-left ones. A few free parameters are introduced to control the displacement of the box. No algorithm is given to determine these parameters. The approach described is only applicable to the minimization problem when only the class function *1S* is involved.

We suggest another strategy to seek the Pareto frontier. For the sake of simplicity, we assume initially that only a minimization problem is considered (all class functions are *1S*) and the problem is solved in the objective space $\mathbf{Y}$. The generalization on the arbitrary case is given subsequently.

Similar to the NC method, let us consider the utopia plane created by anchor points $\boldsymbol{\mu}_i^*$. It is well known that any point $\mathbf{p}$ belonging to the interiority of a convex polygon spanned by $n_{sc}$ vertexes $\boldsymbol{\mu}_i^*$ can be represented as follows:

$$\mathbf{p} = \sum_{i=1}^{n_{sc}} \alpha_i \boldsymbol{\mu}_i^* \quad (24)$$

where the parameters $\alpha_i$ must satisfy the following conditions:

$$0 \leq \alpha_i \leq 1, \quad (25)$$

$$\sum_{j=1}^{n_{sc}} \alpha_j = 1$$

In this approach the notion of the anchor point is used. Although not crucial, an important general remark on the definition of the anchor point should be made. As it was mentioned above, the standard definition assumes an anchor point $\boldsymbol{\mu}_i^*$ corresponds to the solution of the single-optimization problem *min $F_i$* in feasible criterion space $\mathbf{Y}^*$ (see, e.g.,[9]). This definition allows the anchor point corresponding to some objectives to be non-unique. Furthermore, it may not even belong to the boundary of $\mathbf{Y}^*$. Such an example is given below. We suggest the following specification which guarantees the uniqueness of the anchor point for each objective. If the solution of the problem $\boldsymbol{\mu}_i^* = \mathbf{F}(\mathbf{X}^{i*}) \ \{\mathbf{X}^{i*} : \mathbf{X}^{i*} = \arg\min_{\mathbf{Y}^*} F_i\}$ is not unique, then the point corresponding to the minimal values of the other design metrics is to be chosen. It may lead to the problem of trade-off minimization for the remaining objectives. To avoid this, priority in the minimization is introduced. First, $F_i$ is to be minimized, then $F_{i+1}$ and so on up to $F_{i-1}$. The prioritization is introduced in a circular order: *i+1, i+2, …, $n_{sc}$, 1, 2, …, i-1*. A *k-th* prioritization assumes that the *k-th* minimization must not violate all the previous *k-1* ones. At this definition, it is easy to prove that all anchor points belong to the Pareto frontier because they are on the boundary of the feasible space $Y^*$ and no objectives can be improved without deterioration of any other objective.

The free vector $\mathbf{F}^{(0)} = (F_1^{(0)}, F_2^{(0)}, ..., F_{n_{sc}}^{(0)})^T$ is determined in the following way. Let us consider the box *D*. The box *D* is shifted in such a way that its vertex *M* corresponding to the maximal values of the design metrics $(M = (F_{1,\max}, F_{2,\max}, ..., F_{n_{sc},\max})^T)$ lies in the utopia plane (see Figure 2). This means

$$\mathbf{F}^{(0)} = \sum_{i=1}^{n_{sc}} \alpha_i \boldsymbol{\mu}_i^* - \mathbf{a}, \tag{26}$$

where

$$\mathbf{F}^{(0)} = (F_1^{(0)}, F_2^{(0)}, ..., F_{n_{sc}}^{(0)})^T,$$
$$\mathbf{a} = (a_1, a_2, ..., a_{n_{sc}})^T.$$

An even distribution of the coefficients $\alpha_i$ gives us an even distribution of the Pareto set. In comparison to the NC and NBI methods, the approach described below allows us to generate the complete Pareto frontier considering only non-negative coefficients $\alpha_i$ from (25).

One of the possible algorithms for calculation of the coefficients $\alpha_i$ is given in[9] where the following induction procedure is used. First, a uniform distribution of coefficient $\alpha_1$ is considered. The sum of the rest coefficients $\sum_{j=2}^{n_{sc}} \alpha_j$ equals to $1 - \alpha_1$ for each selected value of $\alpha_1$. Then, a uniform distribution of the coefficient $\alpha_2$ is considered for each of these variants and so on until either the last coefficient $\alpha_{n_{sc}}$ is reached or the sum of the coefficients already determined equals to 1. In the latter variant the remaining coefficients equal zero.

It will be shown that the distribution of the Pareto set may be sensitive to the displacement of the box $D$ along the utopia plane especially if the Pareto frontier in concave. To avoid this, generalization of the class functions is performed. It allows us to shrink the search domain substantially.

In order to shrink the search domain defined by the hypercube $H$, it is suggested to introduce generalized class-functions as follows:

$$\overline{F}(\tilde{F}_i) \quad (i = 1, ..., n_{sc}), \tag{27}$$

where $\tilde{F}_i$ is defined by an affine transform

$$\tilde{F}_i = F_j B_{ji}, \quad (i, j = 1, ..., n_{sc}). \tag{28}$$

In the objective space **Y**, it is equivalent to the introduction of a new coordinate system with the basic vectors

$$\mathbf{a}_i = A_{ij} \mathbf{e}_j, \quad (i, j = 1, ..., n_{sc}), \tag{29}$$
$$A^{-1} = B$$

where $\mathbf{e}_j$ ($j = 1, ..., n_{sc}$) are the basic vectors of the original coordinate system.

Then, the search domain can be changed as shown in Figure 3. In particular, it is possible to choose basic vectors $\mathbf{a}_i$ ($i = 1, ..., n_{sc}$) which form an angle $\gamma_c$ to a selected direction **l**. The 2D case is shown in Figure 3. Matrixes **A** and **B** can be easy determined as follows:

$$A = \begin{pmatrix} \cos\gamma_- & \sin\gamma_- \\ \cos\gamma_+ & \sin\gamma_+ \end{pmatrix}, \quad B = \frac{1}{\sin 2\gamma_c} \begin{pmatrix} \sin\gamma_+ & -\sin\gamma_- \\ -\cos\gamma_+ & \cos\gamma_- \end{pmatrix}, \tag{30}$$

where $\gamma_+ = \gamma_n + \gamma_c$, $\gamma_- = \gamma_n - \gamma_c$, $\mathbf{l} = (\cos\gamma_n, \sin\gamma_n)^T$.

In the general case of $R^{n_{sc}}$, it is possible to suggest the following algorithm. There are the following conditions on $\mathbf{a}_i$:

$$(\mathbf{a}_i, \mathbf{l}) = \cos\gamma_c \qquad (i = 1, ..., n_{sc}) \tag{31}$$

All the vectors $\mathbf{a}_i$ are parallel to the lateral area of the hypercone having the angle $\gamma_c$ and the axis along vector **l**. It is important to guarantee a spread distribution of these vectors. At least, the basis created by these vectors must not vanish. It appears possible to obtain even a fully uniform distribution of the basis vectors of such a polyhedral cone.

First, suppose that **l** is directed as follows:

$$\mathbf{l} = \mathbf{l}_0, \tag{32}$$

$$\mathbf{l}_0 = (l_0, l_0, l_0, ..., l_0)^T$$

Assuming that **l** is a unit vector, we obtain its coordinates:

$$l_0 \equiv \cos\gamma_0 = \frac{1}{\sqrt{n_{sc}}} \tag{33}$$

The basis vector $\mathbf{a}_i$ can be determined in the plane created by the vectors $\mathbf{e}_i$ and $\mathbf{l}_0$ (see Figure 4). It is possible to show that

$$\mathbf{a}_i = \frac{\sin\gamma_c}{\sin\gamma_0}\mathbf{e}_i + \frac{\sin(\gamma_0 - \gamma_c)}{\sin\gamma_0}\mathbf{l}_0 \tag{34}$$

From (29), (33) and (34) we obtain

$$A = \frac{\sin\gamma_c}{\sin\gamma_0}I + \frac{\sin(\gamma_0 - \gamma_c)\cos\gamma_o}{\sin\gamma_0}E, \tag{35}$$

where all elements of the matrix $E$ are unities: $\|E_{ij}\| = 1$.

In many cases the shrinking around the lines parallel to the vector $\mathbf{l}_0$ is already sufficient. Nevertheless, it is important to obtain matrix $A$ in the general case of an arbitrary unit vector **l**. For this purpose, it is enough to perform a linear transform mapping the previous pattern in such a way that the vector $\mathbf{l}_0$ is mapped onto the **l**. This purpose is reached by multiplying both parts of equation (34) by an orthogonal matrix $R$:

$$R\mathbf{l}_0 = \mathbf{l} \tag{36}$$

Then, we obtain the basis of vectors $\{\mathbf{a}'_i\}$ ($i = 1, ..., n_{sc}$) uniformly distributed in the lateral area of the hypercone having the axis parallel to the vector **l**:

$$\mathbf{a}'_i = \frac{\sin\gamma_c}{\sin\gamma_0}\mathbf{e}'_i + \frac{\sin(\gamma_0 - \gamma_c)}{\sin\gamma_0}\mathbf{l}, \tag{37}$$

where $\mathbf{e}'_i = R\mathbf{e}_i$ are the components of the Cartesian coordinate system in which the vector **l** has equalled components. It is easy to see that the columns of transition matrix $R$ are the coordinates of the vectors $\mathbf{e}'_i$ in the basis $\{\mathbf{e}_j\}$. Since the transform is orthogonal, all angles are preserved. In particular, $(\mathbf{a}'_i, \mathbf{l}) = \cos\gamma_0$. Now we can write the matrix $A$ in the general form:

$$A = \frac{\sin\gamma_c}{\sin\gamma_0}R^T + \frac{\sin(\gamma_0 - \gamma_c)}{\sin\gamma_0}E, \tag{38}$$

where $\|E_{ij}\| = \|l_j\|$.

If $\gamma_c = \gamma_0$, obviously $\mathbf{a}'_i = \mathbf{e}'_i$ that means the transform becomes orthogonal and is only reduced to a turn of the original Cartesian coordinate system. As a consequence in this case, the matrix $A$ is orthogonal and $B = A^T$.

The general presentation requires the calculation of the orthogonal matrix $R$, the components of which must satisfy the following additional requirements:

$$\cos\gamma_0 \sum_{j=1}^{n_{sc}} R_{ij} = l_i \tag{39}$$

Matrix $R$ is not unique. The simplest way to obtain it is to consider the rotation from the vector $\mathbf{l}_0$ to the vector $\mathbf{l}$ in a Cartesian coordinate system related with these vectors so that

$$R = DT_R D^{-1}, \tag{40}$$

where $T$ is an elementary rotation matrix describing the rotation in the plane created by the first two basis vectors

$$T_R = \begin{bmatrix} (\mathbf{l}_0,\mathbf{l}) & -\sqrt{1-(\mathbf{l}_0,\mathbf{l})} & 0 & \dots & 0 \\ \sqrt{1-(\mathbf{l}_0,\mathbf{l})} & (\mathbf{l}_0,\mathbf{l}) & 0 & \dots & 0 \\ 0 & 0 & 1 & \dots & \dots \\ \dots & \dots & \dots & \dots & \dots \\ 0 & 0 & 0 & \dots & 1 \end{bmatrix}, \tag{41}$$

and the matrix $D$ is the transition matrix from the original basis to some orthogonal basis $\{\mathbf{b}_i\}$: $\mathbf{b}_1 = \mathbf{l}_0$, $\mathbf{b}_2 = \dfrac{\mathbf{l} - (\mathbf{l},\mathbf{l}_0)\mathbf{l}_0}{\sqrt{1-(\mathbf{l},\mathbf{l}_0)^2}}$, ... . The remaining basis vectors $\mathbf{b}_3, \mathbf{b}_4, \dots, \mathbf{b}_{n_{sc}}$ can be easy obtained by the Gram-Schmidt orthogonolization procedure.

Obviously the vectors $\mathbf{a}_i$ ($i = 1, \dots, n_{sc}$) create a basis in $\mathbf{Y}$ which does not vanish. The basis vectors form a search cone similar to the 2D case shown in Figure 3. In the 2D case, $n_{sc} = 2$, $\gamma_0 = \pi/4$ and we obtain formulas (30).

The boundaries of the preference ranges are mapped in according to (28):

$$\tilde{F}_{ik} = F_{jk} B_{ji}, \quad (i, j = 1, \dots, n_{sc}; k = 1, \dots, 5) \tag{42}$$

Transform (27), (28) allows us to shrink the search domain and focus on a much smaller area on the Pareto surface. It makes the algorithm more flexible and much less sensitive to the displacement of box $D$. The transform shrinks the search domain, much like a "light beam" which emits from point $M$ and highlights a spot on the boundary of feasible space $\mathbf{Y}^*$. The direction of the search is easy to handle. It is natural to choose this direction (vector -$\mathbf{l}$) in alignment with the normal to the utopia hyperplane towards the decrement of the objective functions. If no solution is found, the direction is switched on the opposite one. The box $\tilde{D}$, which is the image of the box $D$ after transform (27), (28) is translated along the vector **l** until the whole box $\tilde{D}$ intersects the utopia plane. In this case, the new pseudo-preference $\tilde{F}_{i1}$ corresponds to the former pseudo-preference $\tilde{F}_{i5}$ and the search domain is limited by the box $\tilde{D}$ only.

It should be emphasised that the general representation of matrix $A$ can play a substantial part in seeking the Pareto set nearby its boundary. As it was mentioned above, if we consider orthogonal projection of the Pareto set onto the utopia hyperplane, the images of some Pareto points may not belong to the interiority of a convex polygon spanned by the $n_{sc}$ vertexes $\boldsymbol{\mu}_i^*$ (24), (25). This fact was first noted in[9]. One of the possibilities to resolve this problem, suggested in[9] for the NC method, is based on the use of negative coefficients $\alpha_i$. However, this will cause another problem, this time with the lower evaluation of the coefficients. Inevitably some points in the utopia plane corresponding to negative coefficients may not be orthogonal images of any Pareto point. Another opportunity, offered specifically by our modified PP is described below.

Let us consider the edge vectors of polygon (24), (25): $\mathbf{v}_i = \boldsymbol{\mu}_{i+1} - \boldsymbol{\mu}_i$ ($i = 1, ..., n_{sc}-1$). The point $\mathbf{p_i}$ belongs to a $k$-th edge of the polygon iff $\alpha_m = 0$ ($m \neq k, k+1$). Assuming that vector **l** is related with the normal of the utopia hyperplane. Then, when point $M$ lies on an edge of the polygon, vector **l** is rotated in the direction opposite to the polygon. In another words, **l** is changed in such a way that the orthogonal projection of the end of the vector, drawn from an edge, onto the utopia hyperplane must not fall in the interiority of the polygon. For this purpose, in the utopia plane we introduce a unit vector which is the outer normal to the edge considered. The vector can be defined as:

$$\mathbf{s}_i = \frac{\mathbf{v}_{i-1} + \beta_i \mathbf{v}_i}{|\mathbf{v}_{i-1} + \beta_i \mathbf{v}_i|}, \quad \beta_i = -\frac{(\mathbf{v}_{i-1}, \mathbf{v}_i)}{(\mathbf{v}_i, \mathbf{v}_i)}. \tag{43}$$

Then, the current vector $\mathbf{l}_r$ is determined via $\mathbf{s}_i$ and normal **n** to the utopia hyperplane towards the utopia point $(F_1(\mu_1), F_2(\mu_2), ..., F_{n_{sc}}(\mu_{n_{sc}}))$ as follows:

$$\mathbf{l}_r = -\cos\theta_r \mathbf{n} + \sin\theta_r \mathbf{s}_i, \tag{44}$$
$$0 < \theta_r < \pi/2.$$

The angle $\theta_r$ is a parameter. Changing $\theta_r$ from $\pi/2$ to 0, the vector $\mathbf{l}_r$ is turned from the normal vector $-\mathbf{n}$ to the vector $\mathbf{s}_i$ (see Figures 5a and 5b). Thus, the algorithm is formulated in the following way. If point $M$ in (26) belongs strictly to the interiority of the polygon, then the vector $\mathbf{l}_r$ coincides with the normal $-\mathbf{n}$. If point $M$ lies on an edge of the polygon, then an additional rotation of the vector may be required. To obtain an even distribution of the Pareto set, the number of additional points $N_r$ related with the rotation of a vector $\mathbf{l}_r$ depends on the distance to the vertexes of the edge. For example, the rotation is not required at the anchor points. Generally speaking, it is reasonable to choose the maximal value of $N_r$ at the centre of an edge. The following evaluation of $N_r$ is suggested for a *k-th* edge:

$$N_r = \text{int}(4m\alpha_k \alpha_{k+1}) \quad (m \geq 1) \tag{45}$$

Finally, it is worth noting that this number can be substantially optimized if the information on the current local distribution of the Pareto set is taken into account. For example, if a Pareto solution appears to be at an edge of the polygon, no additional rotation is needed and $N_r = 0$.

The approach described above guarantees catching the entire Pareto surface. In some very special cases, like the one shown in Figure 6, local Pareto points (e.g., a point **P**) may be obtained which are not global Pareto solutions. The algorithm excludes such points easily. Let us locate the box *D* in such a way that the point *M* is at some point **P** investigated as a candidate Pareto solution and set *A = I* in (29), (35). If the point is a global Pareto solution (e.g., a point **P′**), no any other solution can be obtained. It immediately follows from the contact theorem[21]. Thus, we have a criterion for verification if the solution is a global Pareto solution.

To avoid undesirable severe skewing of the search domain in the algorithm one may recommend preliminary scaling of the objective functions similar to[1]:

$$F_i^{sc} = \frac{F_i - F_{i,\min}}{F_{i,\max} - F_{i,\min}} \tag{46}$$

**IV. Generalized Pareto set**

Now we propose the generalization of the approach described above in the arbitrary case when all class functions *1S-4S* are involved.

First, let us introduce the notion of reference class-functions. A class-function $\overline{F}(F_i)$ is said to be a reference class-function if its preferences include the entire interval of the variation of the objective $F_i$ in the space $\mathbf{Y}^*$. We map the objective space $Y^*$ onto the class-function space $\mathbf{Z}^*$ by the reference class-functions $\overline{F}(F_i)$ (*i = 1, …, $n_{sc}$*): $Y \mapsto Z$ ($\mathbf{Z} \in R^{n_{sc}}$). The solution belonging to a Pareto set in the space $\mathbf{Z}^*$ can be interpreted as a generalized Pareto optimal (GPO) solution[17]. Vector $\mathbf{a} \in \mathbf{X}^*$ is said to be GPO solution iff it does not exist any vector $\mathbf{b} \in \mathbf{X}^*$ such that $\overline{F}(F_i(\mathbf{b})) \leq \overline{F}(F_i(\mathbf{a}))$ $(i=1,...,n_{sc})$ and there is $l \leq M$: $\overline{F}(F_l(\mathbf{b})) < \overline{F}(F_l(\mathbf{a}))$ where $\overline{F}(F_i)$ are the reference class-functions. The generalization is naturally related with the PP formulation of the optimization problem based on the introduction of such nonmonotonic classes as *3S* and *4S*. In[17], it was proved that any solution of optimization problem (5)-(7) is a generalized Pareto solution. If the problem is considered in the standard form for minimization (1)-(3), then it is easy to prove the following statement.

*Proposition 2*: Assume that the minimization problem (1)-(3) is solved. If the vector **a** is a Pareto optimal solution then it is a GPO solution and vice a versa.

*Proof*: Suppose that **a** is not GPO solution, then there is a vector **b** such that $\overline{F}(f_i(\mathbf{b})) - \overline{F}(f_i(\mathbf{a})) \leq 0$ $(i = 1, ..., M)$ and there is $1 \leq M$: $\overline{F}(f_l(\mathbf{b})) - \overline{F}(f_l(\mathbf{a})) < 0$. Meanwhile, since the **a** is a Pareto optimal solution there is $k \leq M$: $f_k(\mathbf{a}) - f_k(\mathbf{b}) < 0$. From the mean value theorem there is $C > 0$: $\overline{F}(f_k(\mathbf{a})) - \overline{F}(f_k(\mathbf{b})) = C(f_k(\mathbf{a}) - f_k(\mathbf{b})) < 0$. It follows immediately from the fact that the class function is smooth and belongs to *1S*. Thus, we reach a contradiction. The proof of the converse statement is similar because the inverse function $\overline{F}^{-1}$ is also smooth and its derivative is strictly positive. □

Thus, a classical Pareto solution is always a GPO solution. The problem of seeking GPO solutions is associated only with minimization. The appropriate class-functions can be considered as new objectives and then the problem can be solved by the algorithm described above. Following this approach, only the *1S* class-function is to be introduced for the new objectives which are the original class-functions.

Due to the different nature of the optimization problems associated with classes, *2S-4S*, the definition of the anchor points cannot be based on minimization of an objective. Instead, we shall consider the minimization of a reference class function:

$$\overline{\boldsymbol{\mu}}_i^* = \overline{\mathbf{F}}(\mathbf{F}(\mathbf{X}^{i*})) \ \{\mathbf{X}^{i*} : \mathbf{X}^{i*} = \arg\min_{\mathbf{z}^*} \overline{F}(F_i(\mathbf{X}))\} \tag{47}$$

Without loss of generality, assume that point $F_{i1}$ corresponds to the appropriate anchor point $\boldsymbol{\mu}_i^* = \mathbf{F}(\mathbf{X}^*)$ of the *i-th* soft metric for all reference class-functions apart from *4S*. For reference class-function *3S* this correspondence follows from its definition. In the case of class-function *4S*, the anchor point is assumed to be at the middle of interval $[F_{1L} F_{1R}]$.

Now we can consider a tradeoff matrix $\overline{T}$ as follows:

$$T = \begin{bmatrix} \overline{F}_1 & \overline{F}(F_{12}) & ... & \overline{F}(F_{1n_{sc}}) \\ \overline{F}(F_{21}) & \overline{F}_1 & ... & \overline{F}(F_{2n_{sc}}) \\ ... & ... & ... & ... \\ \overline{F}(F_{n_{sc}1}) & \overline{F}(F_{n_{sc}2}) & ... & \overline{F}_1 \end{bmatrix}, \tag{48}$$

In matrix $\overline{T}$ an $i$-th row represents the coordinates of the anchor point $\overline{\boldsymbol{\mu}}_i^*$. In the main diagonal all elements are equal to $\overline{F}_1$ because of the choice of the preferences $F_{i1}$. Here, it is assumed that for the 4S class-function the difference between $\min \overline{F}(F)$ and $\overline{F}(F_1)$ can be neglected.

Then, the pseudo nadir points in the space $\mathbf{Z}$ can be defined as $\overline{F}_{i,\max} = \max_j \overline{F}(F_{ij})$. Obviously, $\overline{F}_{i,\max} \leq \overline{F}_5$. A hypercube $\overline{H}$ limiting the investigation of feasible space $\mathbf{Z}^*$ is defined as follows: $\overline{H} = [\overline{F}_1 \overline{F}_{1,\max}] \times [\overline{F}_1 \overline{F}_{2,\max}] \times ... \times [\overline{F}_1 \overline{F}_{n_{sc},\max}]$.

If a class-function is determined, it maps the pseudo-preference vector $\mathbf{P}_i = (F_{i1}, F_{i2}, F_{i3}, F_{i4}, F_{i5})^T$ onto the vector $\overline{\mathbf{P}}_i = (\overline{F}_{i1}, \overline{F}_{i2}, \overline{F}_{i3}, \overline{F}_{i4}, \overline{F}_{i5})^T$, where $\overline{F}_{ip} = \overline{F}(F_{ip})$ ($p = 1, 2, ..., 5$). The inverse function $\overline{F}^{-1}(F_i)$ maps the vector $\overline{\mathbf{P}}_i$ onto the pseudo-preference vector $\mathbf{P}_i = (\overline{F}^{-1}(\overline{F}_{i1}), \overline{F}^{-1}(\overline{F}_{i2}), \overline{F}^{-1}(\overline{F}_{i3}), \overline{F}^{-1}(\overline{F}_{i4}), \overline{F}^{-1}(\overline{F}_{i5}))^T$ which is not unique for the classes *3S* and *4S*. Anyway, the determination of the vector $\overline{\mathbf{P}}_i$ ($i = 1, ..., n_{sc}$) leads to the determination of the pseudo-preference vector $\mathbf{P}_i$.

Let us determine the vector $\overline{\mathbf{P}}_i$ as follows:

$$\overline{\mathbf{P}}_i \equiv (\overline{F}_{i1}, \overline{F}_{i2}, \overline{F}_{i3}, \overline{F}_{i4}, \overline{F}_{i5})^T = \overline{F}_i^{(0)}(1,1,1,1,1)^T + \overline{a}_i(0, \frac{1}{4}, \frac{1}{2}, \frac{3}{4}, 1)^T, \tag{49}$$

where $\overline{a}_i = \overline{F}_{i5} - \overline{F}_1$.

A free vector $\overline{\mathbf{F}}_i^{(0)}$ ($i = 1, ..., n_{sc}$) is determined in the following way. Assuming that in the space $\mathbf{Z}$ a box $\overline{D}$ is the analogue of the box $D$ in $\mathbf{Y}$. The box $\overline{D}$ is shifted in such a way that its vertex $\overline{M}$ corresponding to $(\overline{F}_{15}, \overline{F}_{25}, ..., \overline{F}_{n_{sc}5})$ lies on the utopia plane in $\mathbf{Z}$ (see Figure 7a). This means

$$\overline{\mathbf{F}}^{(0)} = \sum_{i=1}^{n_{sc}} \alpha_i \overline{\boldsymbol{\mu}}_i^* - \overline{\mathbf{a}}, \tag{50}$$

where

$$\overline{\mathbf{F}}^{(0)} = (\overline{F}_1^{(0)}, \overline{F}_2^{(0)}, ..., \overline{F}_{n_{sc}}^{(0)})^T,$$

$$\overline{\mathbf{a}} = (\overline{a}_1, \overline{a}_2, ..., \overline{a}_{n_{sc}})^T.$$

If a class-function $\overline{F}(F_i)$ is determined, it maps the pseudo-preference vector $\mathbf{P}_i = (F_{i1}, F_{i2}, F_{i3}, F_{i4}, F_{i5})^T$ onto the vector $\overline{\mathbf{P}}_i = (\overline{F}_{i1}, \overline{F}_{i2}, \overline{F}_{i3}, \overline{F}_{i4}, \overline{F}_{i5})^T$ where $\overline{F}_{ip} = \overline{F}(F_{ip})$ ($p = 1, 2, ..., 5$). The inverse function $\overline{F}^{-1}(F_i)$ maps the vector $\overline{\mathbf{P}}_i$ onto the pseudo-preference vector $\mathbf{P}_i = (\overline{F}^{-1}(\overline{F}_{i1}), \overline{F}^{-1}(\overline{F}_{i2}), \overline{F}^{-1}(\overline{F}_{i3}), \overline{F}^{-1}(\overline{F}_{i4}), \overline{F}^{-1}(\overline{F}_{i5}))^T$ which is not unique for the classes *3S* and *4S*. In any case, the determination of the vector $\overline{\mathbf{P}}_i$ ($i = 1, ..., n_{sc}$) leads to the determination of the pseudo-preference vector $\mathbf{P}_i$.

Thus, moving the box $\overline{D}$ in $\mathbf{Z}$ we automatically define the pseudo-preferences for the original class functions. For example, in the 2D case, if box $\overline{D}$ is in the lower-right angle of the hypercube $\overline{H}$ then the metric $F_2$ is optimized as a trade-off of the metric $F_1$, and it does not matter which class-functions are involved. It leads to the full possible range of $F_1$ in the definition of the pseudo-preference vector $P_1$ from $F_{1,min}$ to $F_{1,max}$ while the metric $F_2$ can only be changed in the vicinity of the anchor point. This is illustrated in Figure 7b. For the sake of determination, it is assumed that $\overline{F}_2$ is the 3S class-function.

### V. Test cases

The method described above is validated using a few test cases. The test cases include examples with both Pareto convex and concave frontiers. It is shown that the standard realization of the PP may lead to a high sensitivity of the location of the Pareto points to the displacement of box $D$ (or $D^*$) in the case of a concave Pareto frontier.

**Example 1:**

First, we consider the following simple algebraic test cases:

$$\min (x, y)^T \qquad (51)$$

case one, constraints:

$$x^2 + y^2 \leq 1, \qquad (52)$$
$$x > -1, y > -1$$

case two, constraints:

$$x^2 + y^2 \geq 1, \qquad (53)$$
$$x \geq 0, y \geq 0$$

In the first case, (51) and (52), the Pareto surface is convex while in the second case, (51) and (53), it is concave. The latter case also represents an example where the standard definition of an anchor point leads to non-uniqueness. In this case any point with coordinates $x = 0, y \geq 1$ is a solution of the single-objective minimization problem: $\min_{x,y} x$. In our formulation we have only two anchor points in total: (0,1) and (1,0).

The solution of problem (51), (52) is shown in Figure 8a. An even displacement of point *M* (or box *D*) along the utopia line does not lead to a completely even distribution of the Pareto points since there are gaps nearby the anchor points. In the case of shrinking the search domain via transform (27) we obtain an even representation of the Pareto frontier (Figure 8b).

In problem (51), (53) the Pareto frontier is concave. The displacement along the utopia line does not provide any solution, except the utopia points because box *D* is in the unfeasible space $\mathbf{Y}\backslash\mathbf{Y}^*$. If the displacement of point *M* is performed along the line parallel to the utopia line as shown in Figure 9, a complete Pareto frontier can be obtained. Yet, this Pareto set is not evenly distributed. The performance of the algorithm is demonstrated in Figure 10. In this figure, the contour plots of the AOF are located inside the search domain. The surface shows the values of the AOF inside box $\tilde{D}$. The complete solution is represented in Figure 10 along with the contour lines located in the corresponding boxes $\tilde{D}$ for each of the utopia plane points. If a solution is not found (as shown in Figure 10), the direction of the search is switched to the opposite one.

**Example 2:**

This example is similar to the previous test case however, the different constraint set provides an opportunity to demonstrate the optimization in the class-function space **Z**:

$$\min x, \max y \tag{54}$$

under constraints :

$$x^2 + y^2 \leq 1 \tag{55}$$

$$x < 0, y > 0$$

This example is also used to demonstrate the capability of the generalized algorithm based on the reference class-functions. It is worth noting that this problem can be trivially reduced to the minimization problem. Nevertheless, for the purposes of this example, it is deliberately considered as the problem with different categories of *1S* and *2S* types. The obtained solutions are represented in Figure 11.

The next test case is a *3D* one and leads as to an example where the orthogonal images of Pareto point onto the utopia surface are not necessary in the interior of the polygon related with the anchor points.

**Example 3:**

A 3D minimization problem:

$$min \ (x, y, z)^T \qquad (56)$$

under constraints

$$x^2 + y^2 + z^2 \geq 1 \qquad (57)$$

$$x > 0, \ y > 0, \ z > 0$$

The standard definition of the anchor points again leads to non-uniqueness. It is easy to see that in our formulation we have only three anchor points: (0,0,1), (1,0,0) and (0,1,0). Not the entire orthogonal projection of the Pareto surface onto the utopia surface appears to be in the triangle created by the anchor points (see Figure 12). For this reason, a method such as NBI, for example, is not able to catch the entire Pareto frontier if the coefficients $\alpha_i$ in (25) are not negative. Algorithm (44), (45) is used to provide the complete representation of the Pareto frontier. The utopia plane and points, distributed in the polygon (triangle) according to algorithm (24), (25), are shown in Figure 13a. The solution providing the complete Pareto frontier is given in Figure 13b.

**VI. Conclusions**

The Physical Programming method has been modified to generate an even distribution of the entire Pareto set. The modification is based on the generalization of the class-functions which leads to shrinking of the search domain. The orientation of the search domain in space can be easily conducted. It allows the method to provide an even distribution of the entire Pareto surface. The generation is performed for both convex and non-convex Pareto frontiers. The method does not generate non-Pareto solutions. A simple algorithm has been proposed to remove local Pareto solutions which are not the global Pareto solutions. A new, more compact representation of the class-functions has been suggested.

The algorithm has been also extended to the case of the generalized Pareto solutions.

The suggested approach has been verified by different test cases, including the generation of both convex and concave Pareto frontiers.

**Acknowledgements**

This research was founded by European Aerospace Project VIVACE (Value Improvement through a Virtual Aeronautical Collaborative Enterprise).

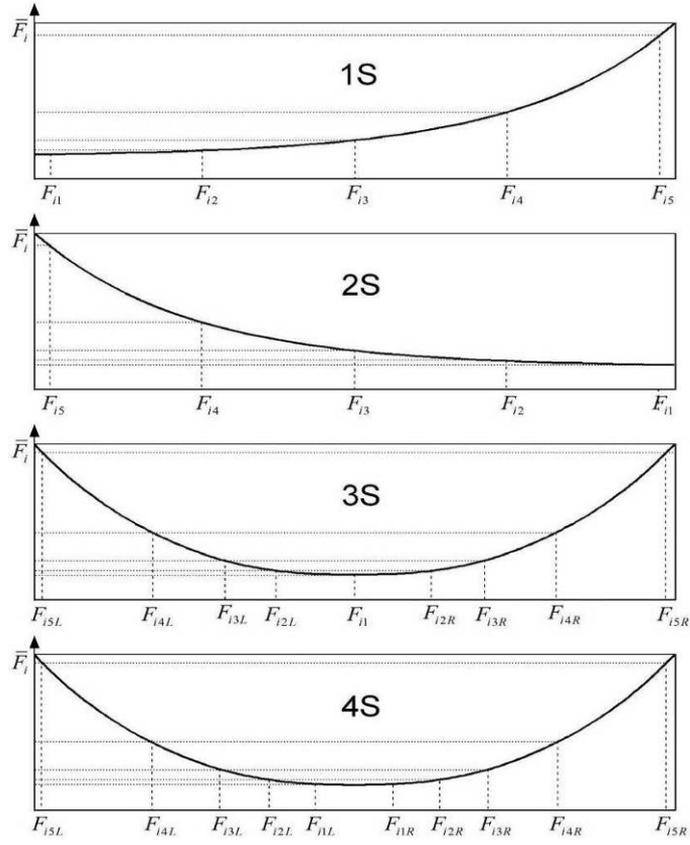

Fig. 1. The class-functions

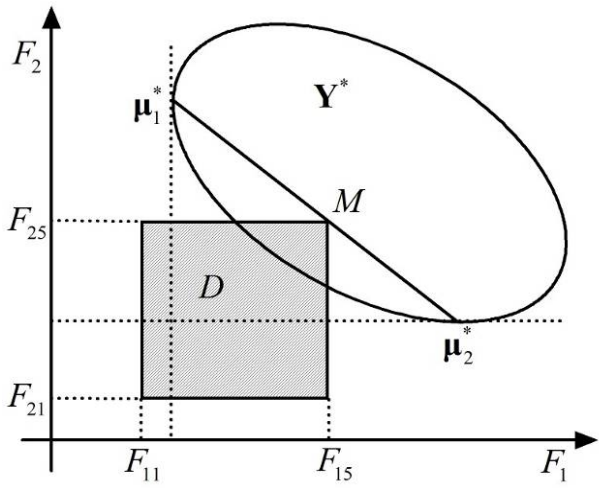

Fig. 2. Original search domain in the objective space

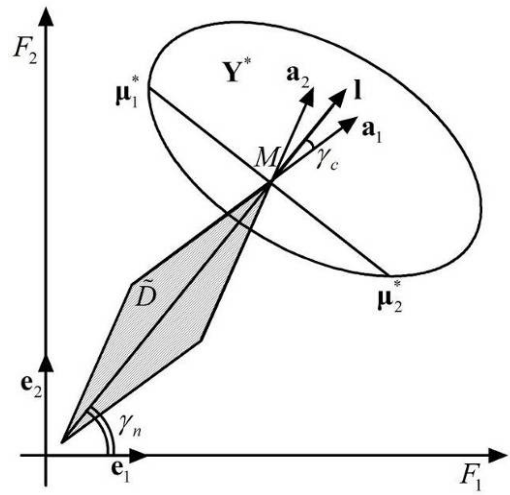

Fig. 3. New search domain after transform

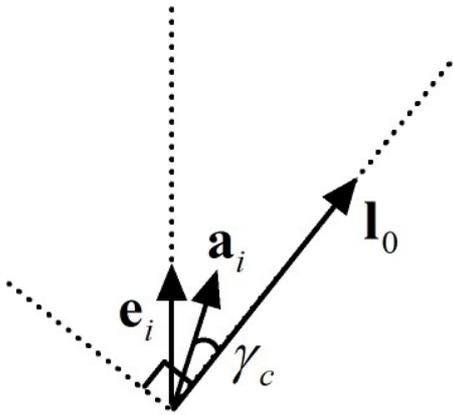

Fig. 4. Local basis vectors in an i-th hyperplane

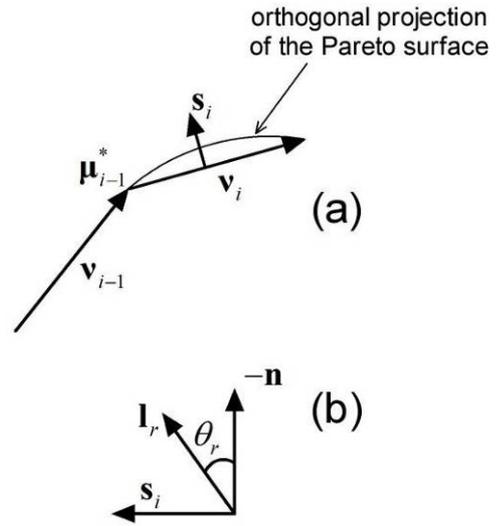

Fig. 5. Rotation of the search domain

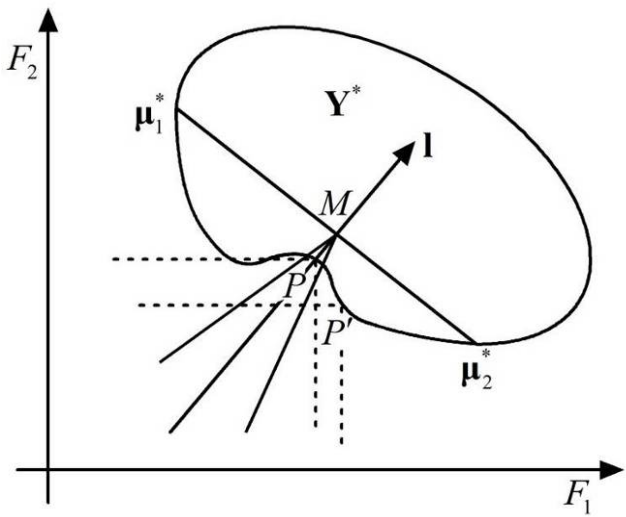

Fig. 6. Verification of global Pareto

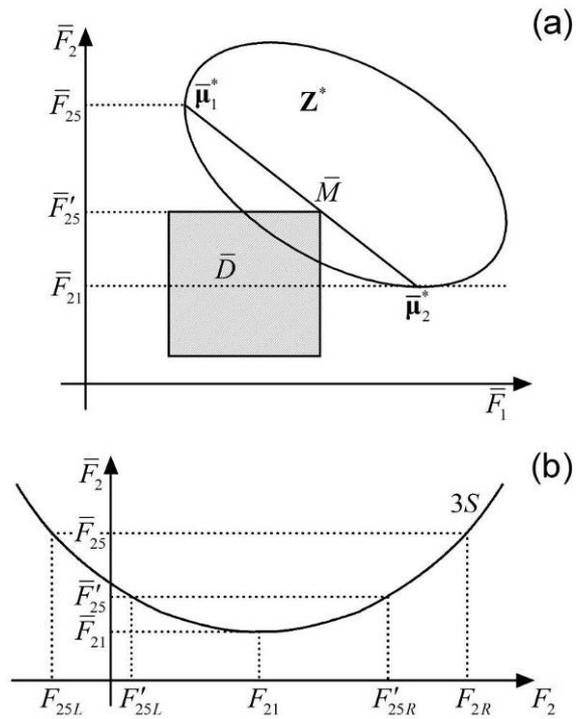

Fig. 7. Search domain in the class-function space and its effect on pseudo-preferences

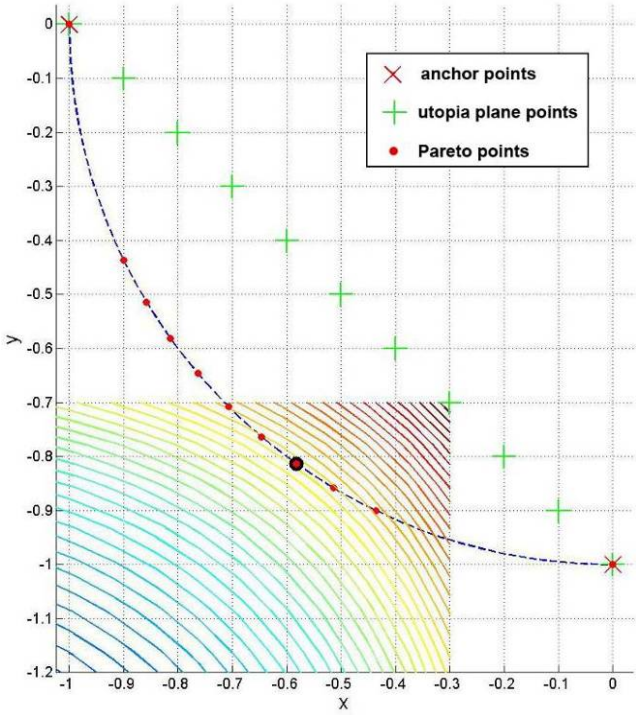

Fig. 8a. Standard search domain. Convex Pareto frontier

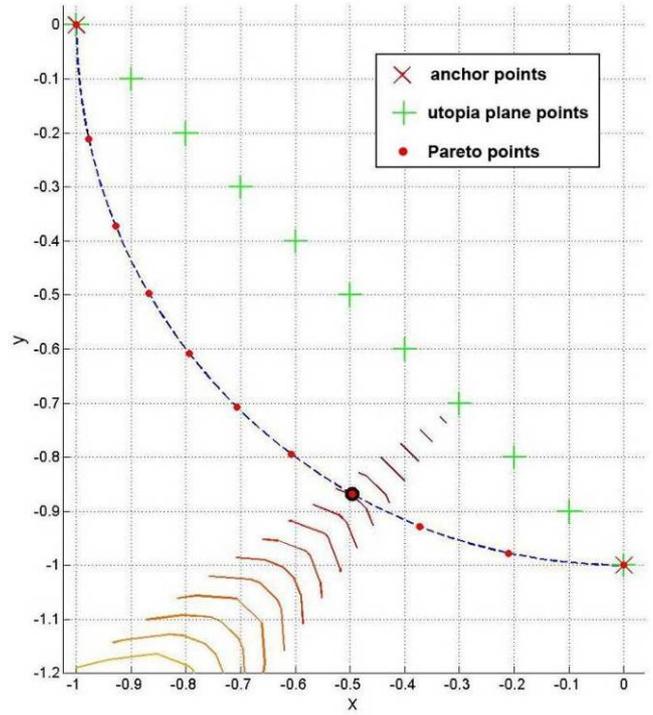

Fig. 8b. Transformed search domain. Convex Pareto frontier

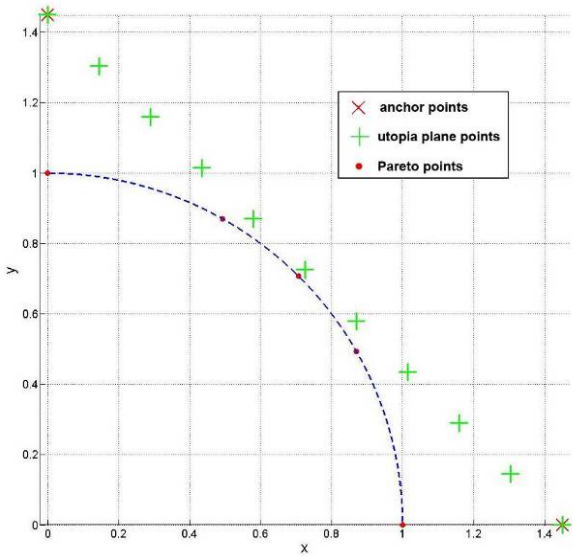

Fig. 9. Standard search domain. Concave Pareto frontier

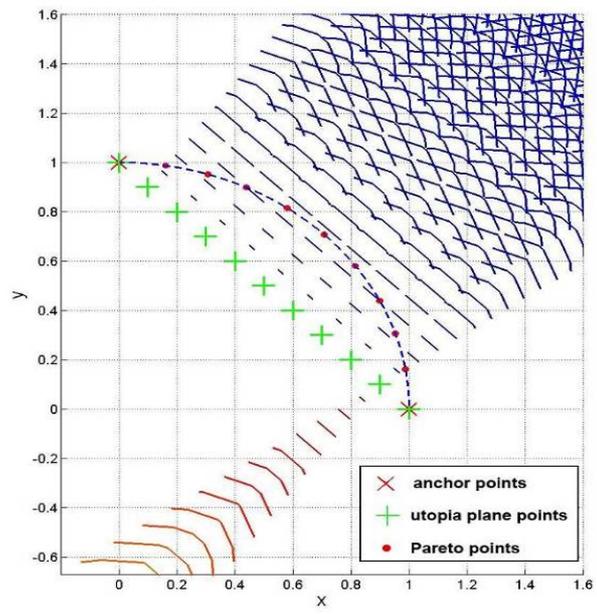

Fig. 10. Transformed search domain. Concave Pareto frontier

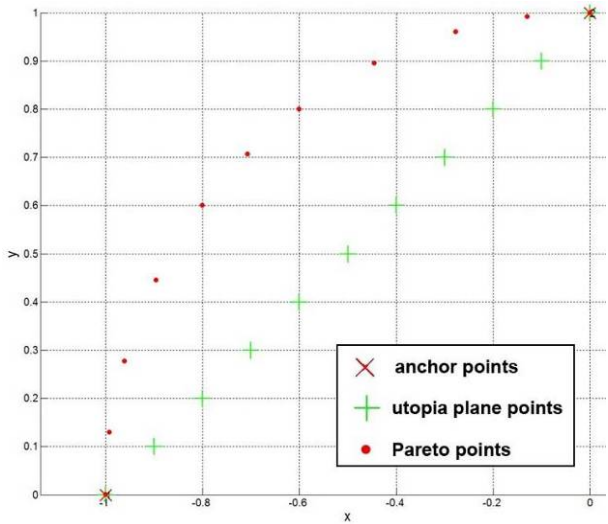

Fig. 11. Pareto frontier. Solving min-max problem using the reference class-function algorithm

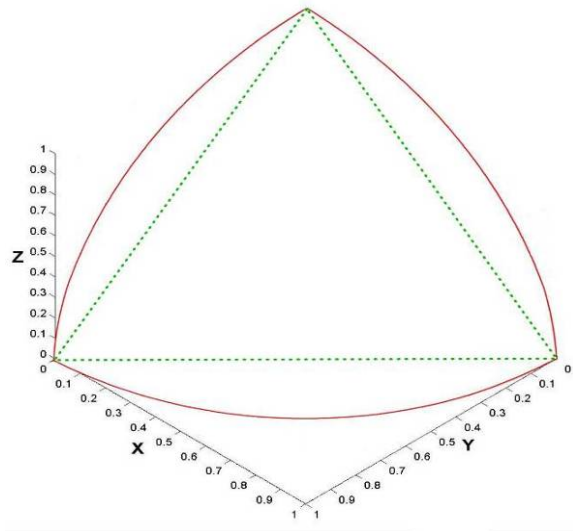

Fig. 12. 3D test case. Projection of the Pareto surface onto the utopia plane

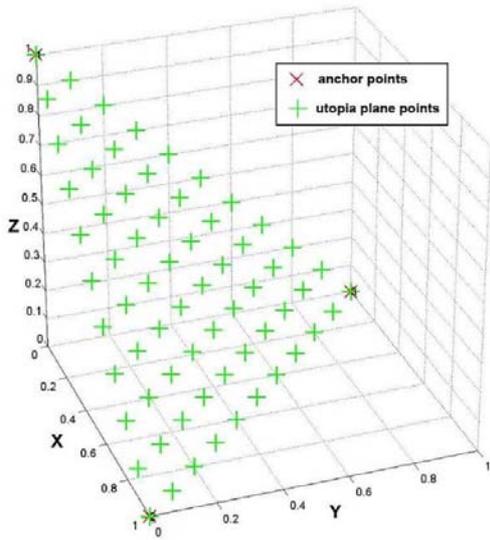

Fig. 13a. 3D test case. Utopia plane

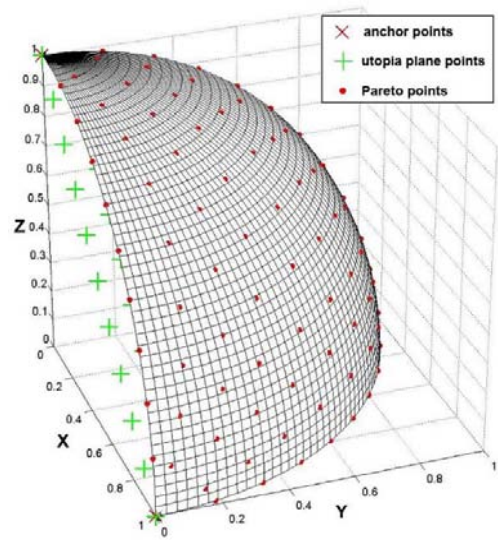

Fig. 13b. 3D test case. Pareto surface